\documentclass[10pt,a4paper]{article}
\usepackage[latin1,applemac]{inputenc}
\usepackage{amsmath, amssymb, amsfonts, amsthm, amscd, epsfig, multicol, marvosym, float,mathrsfs}
\usepackage{mathpazo} 
\usepackage[scaled=.95]{helvet} 
\usepackage{courier}

\usepackage{multicol}

\usepackage[left=2cm,top=3cm,bottom=3cm,right=2cm]{geometry}
\usepackage{fancyheadings}
\usepackage{soul}
\usepackage[usenames]{color}

\begin{document}
\author{Daniele Ritelli \footnote{Dipartimento di Scienze Statistiche, Università di Bologna \textcolor{blue}{{\tt daniele.ritelli@unibo.it}}}}
\title{Another proof of $\displaystyle\boldsymbol{\zeta(2)=\frac{\pi^2}{6}}$ using double integrals}
\date{}
\maketitle

\begin{abstract}
Starting from the double integral
\begin{equation*}
\int_0^{\infty}\int_0^\infty\frac{{\rm d}x{\rm d}y}{(1+y)(1+x^2y)}
\end{equation*}
we give another solution to the Basel Problem
\[
\zeta(2)=\sum_{n=1}^\infty\frac{1}{n^2}=\frac{\pi^2}{6}.
\]

2000 AMS Classification {\bf 40A25}

\vspace{.5cm}


\vspace{.5cm} 

{\sc Note:} To appear on American Mathematical Monthly

\end{abstract}

The celebrated Euler identity, known as the Basel Problem,
\begin{equation}\label{euler}
\sum_{n=1}^\infty\frac{1}{n^2}=\frac{\pi^2}{6}
\end{equation}
has been proved in many different ways.  In this note, we focus on the derivation of \eqref{euler}, taking advantage of the nice interplay between a double integral and a geometric series as has appeared in several articles on this subject, \cite{A, H, Bb}. As Sir
Michael Atiyah declared in an interview \cite{[1]}:
\begin{quote}
\textit{Any good theorem should have several proofs, the more the better. For
two reasons: usually, different proofs have different strengths and
weaknesses, and they generalize in different directions: they are not just
repetitions of each other}
\end{quote}
we find that there is always something worthy of attention in a new proof of a known result. Here, we provide a proof that uses a rational function with the lowest degree among the functions used in different proofs of the same kind.

The author who inaugurated this approach was Apostol, \cite{A}, inspired by Beukers' paper \cite{B}, where the double integral
\begin{equation}\label{box}
\int_0^1\int_0^1\frac{{\rm d}x{\rm d}y}{1-xy} 
\end{equation}
was utilized to prove the irrationality of $\zeta(2).$ Apostol, instead, evaluated \eqref{box} in two ways. First, by expanding $1/(1-xy)$ into a geometric series and then, with a change of variables corresponding to the rotation of the coordinate axes through the angle $\pi/4$~radians. By equating the expression so obtained, the value of $\zeta(2)$ is found. It is worth noting that for the second evaluation, one has to compute the elaborate integrals
\[
\begin{split}
I_1&=\int_0^{1/\sqrt2}\frac{1}{\sqrt{2-u^2}}\arctan\left(\frac{u}{\sqrt{2-u^2}}\right){\rm d}u,\,\text{and}\\
I_2&=\int_{1/\sqrt2}^{\sqrt2}\frac{1}{\sqrt{2-u^2}}\arctan\left(\frac{\sqrt2-u}{\sqrt{2-u^2}}\right){\rm d}u
\end{split}
\]
using suitable trigonometric changes of variables. The evaluation of Beukers, Calabi and Kolk, \cite{Bb}, is similar. They expand $1/(1-x^2y^2)$ into a geometric series to obtain the $\zeta(2)$ series, after which they introduce the two dimensional trigonometric changes of variables
\[
x=\frac{\sin u}{\cos v},\quad y=\frac{\sin v}{\cos u}
\]
to evaluate the double integral. The proof of Hirschhon \cite{Hi} stems from the double inequality
\begin{equation}\label{hirschhon}
2\left(\arctan\frac{a}{1+\sqrt{1-a^2}}\right)^2<\sum_{n=0}^\infty\frac{a^{4n+2}}{(2n+1)^2}<2\left(\arctan a\right)^2,\,0<a<1,
\end{equation}
and a passage to the limit as $a\to1.$ The derivation of \eqref{hirschhon} is based on an integral inequality with regard again to the function $f(x,y)=1/(1-x^2y^2).$ To obtain \eqref{hirschhon}, two integrals of $f(x,y)$ over two different regions of the plane are computed. All these approaches have in common the need to remove a singularity at the point $(1,1)$ of the integrand. Our proof is inspired by \cite{H}, where another definite integral
\begin{equation}\label{hh}
\int_0^{+\infty}\int_0^1\frac{x}{(1+x^2)(1+x^2y^2)}{\rm d}x{\rm d}y
\end{equation} 
is computed, first by integrating with respect to $x$ and then with respect to $y$ and vice versa.  The same integral is considered in the probabilistic proof given in \cite{P}, where integral \eqref{hh} comes from the product of two positive Cauchy random variables. The proof of \cite{H}, as well our proof, uses functions with no singularity in the domain of integration, so we can consider that in same sense these proofs are simpler. Moreover, our proof uses a lower degree rational function than the one used in \cite{H}.

Our starting point, as with most of the papers on this subject, is that \eqref{euler} is equivalent to
\begin{equation}\label{harper}
\sum_{n=0}^\infty\frac{1}{(2n+1)^2}=\frac{\pi^2}{8}.
\end{equation}
In our proof we will show \eqref{harper} starting from the double integral
\begin{equation}\label{george1}
\int_0^{\infty}\int_0^\infty\frac{{\rm d}x{\rm d}y}{(1+y)(1+x^2y)}.
\end{equation}
If we integrate \eqref{george1} first with respect to $x$ and then to $y,$ we find that:
\begin{equation}\label{dueb}
\begin{split}
\int_0^\infty\left(\frac{1}{1+y}\int_0^\infty\frac{{\rm d}x}{1+x^2y}\right){\rm d}y&=\int_0^\infty\left(\frac{1}{1+y}\left[\frac{\arctan(\sqrt{y}\,x)}{\sqrt{y}}\right]_{x=0}^{x=\infty}\right){\rm d}y\\
&=\frac\pi2\int_0^\infty\frac{{\rm d}y}{\sqrt{y}(1+y)}=\frac\pi2\int_0^\infty\frac{2u}{u(1+u^2)}{\rm d}u=\frac{\pi^2}{2}
\end{split}
\end{equation}
where we used the change of variable $y=u^2$ in the last step. 
Reversing the order of integration yields
\begin{equation}\label{duec}
\begin{split}
\int_0^{\infty}\left(\int_0^\infty\frac{{\rm d}y}{(1+y)(1+x^2y)}\right){\rm d}x&=\int_0^\infty\frac{1}{1-x^2}\left(\int_0^\infty\left(\frac{1}{1+y}-\frac{x^2}{1+x^2y}\right){\rm d}y\right){\rm d}x\\
&=\int_0^\infty\frac{1}{1-x^2}\ln\frac{1}{x^2}{\rm d}x=2\int_0^\infty\frac{\ln x}{x^2-1}{\rm d}x.
\end{split}
\end{equation}
Hence, equating \eqref{dueb} and \eqref{duec} we get
\begin{equation}\label{due}
\int_0^\infty\frac{\ln x}{x^2-1}{\rm d}x=\frac{\pi^2}{4}.
\end{equation}
Now split the integration domain in \eqref{due} between $[0,1]$ and $[1,\infty)$ and change the variable $x=1/u$ in the second integral, so that
\begin{equation}\label{tre}
\begin{split}
\int_0^\infty\frac{\ln x}{x^2-1}{\rm d}x&=\int_0^1\frac{\ln x}{x^2-1}{\rm d}x+\int_1^\infty\frac{\ln x}{x^2-1}{\rm d}x\\
&=\int_0^1\frac{\ln x}{x^2-1}{\rm d}x+\int_0^1\frac{\ln u}{u^2-1}{\rm d}u.
\end{split}
\end{equation}
From \eqref{due} and \eqref{tre} we get
\begin{equation}\label{qua}
\int_0^1\frac{\ln x}{x^2-1}{\rm d}x=\frac{\pi^2}{8}.
\end{equation}
Equation \eqref{harper} now follows, expanding, as in \cite{H}, the denominator of the integrand on the left hand side of \eqref{qua} into a geometric series and using the Monotone Convergence Theorem (see \cite{CK} pp. 95-96). Thus, we have: 
\begin{equation}\label{cin}
\int_0^1\frac{\ln x}{x^2-1}{\rm d}x=\int_0^1\frac{-\ln x}{1-x^2}{\rm d}x=\sum_{n=0}^{+\infty}\int_0^1(-x^{2n}\ln x)\,{\rm d}x.
\end{equation}
Integrating by parts yields
\begin{equation}\label{cinq}
\int_0^1(-x^{2n}\ln x)\,{\rm d}x=\left[-\frac{x^{2n+1}}{2n+1}\ln x\right]_0^1+\int_0^1\frac{x^{2n}}{2n+1}{\rm d}x=\frac{1}{(2n+1)^2}
\end{equation}
so that considering \eqref{cinq}, we can write \eqref{cin} as
\begin{equation}\label{cinu}
\int_0^1\frac{\ln x}{x^2-1}{\rm d}x=\sum_{n=0}^{+\infty}\frac{1}{(2n+1)^2}
\end{equation}
and we are done equating \eqref{qua} and \eqref{cinu}.

\bigskip

\noindent\textit{Dipartimento di Statistica, Universit\`{a} di Bologna,\\
Viale Filopanti, 5 Bologna\\
daniele.ritelli@unibo.it}

\end{document}